
\input amstex
\documentstyle{amsppt}
\NoBlackBoxes
\NoRunningHeads
\topmatter
\title
Infinite-Dimensional Geometry of the Universal
Deformation of the Complex Disk \endtitle
\author
Denis~V.~Yur'ev
\endauthor
\address
Institute for Information Technologies,
23 Avtozavodskaya, Moscow 109280, Russia
\endaddress
\date
July 11, 1992
\enddate
\rightheadtext{Infinite-Dimensional Geometry of the Universal
Deformation of the Complex Disk }
\abstract
The universal deformation of the complex disk is studied from the viewpoint of
infinite-dimensional geometry. The structure of a subsymmetric space on the
universal deformation is described. The foliation of the universal deformation
by subsymmetry mirrors is shown to determine a real polarization. \endabstract
\endtopmatter

\document
\eightpoint

The subject of this paper may be of interest to specialists in algebraic
geometry and representation theory as well as to researchers dealing with
mathematical problems of modern quantum field theory.

The universal deformation of the complex disk is one of the crucial concepts
used in the geometric statement of quantum conformal field theory \cite
{1} and quantum-field theory of strings \cite {2} (see also \cite {3}).
The characteristic feature of the approach developed in the present paper is
that the universal deformation of the complex disk is studied in terms of
infinite-dimensional geometry. On this way the structure of the subsymmetric
space \cite {4--6} on the universal deformation is described. The foliation of
the universal deformation defined by the mirrors of subsymmetries determines a
real polarization. For a long time real polarizations on complex manifolds and
their quantization have been attracting the attention of
mathematicians dealing with algebraic geometry and representation theory and of
specialists in mathematical physics \cite {44--47}. The results of this paper
confirm the importance of studying such polarizations
and expose a connection between the traditions of classical synthetic geometry
and recent trends in algebraic geometry, representation theory,
and modern quantum field theory.

\head 1. The infinite-dimensional geometry of the flag manifold of the
Virasoro-Bott group (the base of the universal deformation of the
complex disk).
\endhead

\subhead 1.1. The Virasoro algebra, the Virasoro-Bott group, and the Neretin
semigroup \endsubhead
Let $\operatorname {Diff} (S^1)$ denote the group of diffeomorphisms of the
unit circle $S^1$. The group manifold $\operatorname {Diff} (S^1)$ splits into
two connected components, the subgroup $\operatorname {Diff}_+ (S^1)$ and the
coset $\operatorname {Diff}_- (S^1)$. The diffeomorphisms in
$\operatorname {Diff}_+ (S^1)$ preserve the orientation on the circle $S^1$ and
those in $\operatorname {Diff}_- (S^1)$ reverse it.

The Lie algebra of $\operatorname {Diff}_+ (S^1)$ can be identified
with the linear space $\operatorname {Vect} (S^1)$ of smooth vector fields on
the circle equipped with the commutator
$$
[v(t)d/dt,u(t)d/dt]=(v(t)u'(t)-v'(t)u(t))d/dt.
\tag 1
$$
In the basis
$$
s_n=\sin(nt)d/dt, \quad c_n=\cos(nt)d/dt,
\quad h=d/dt
\tag 2
$$
the commutation relations have the form
$$\aligned
[s_n,s_m] &= 0.5 ((m-n)s_{n+m} + \operatorname{sgn}(n-m)(n+m)s_{|n-m|})\\
[c_n,c_m] &= 0.5 ((n-m)s_{n+m} + \operatorname{sgn}(n-m)(n+m)s_{|n-m|})\\
[s_n,c_m] &= 0.5 ((m-n)c_{m+n} - (m+n)c_{|n-m|})-n\delta_{nm}h\\
[h,s_n] &= nc_n\\
[h,c_n] &= ns_n
\endaligned
\tag 3
$$
The complexification of the Lie algebra $\operatorname {Vect} (S^1)$ will
be denoted by $\operatorname {Vect}_{\Bbb C} (S^1)$. It is convenient to
choose the following basis in $\operatorname {Vect}_{\Bbb C} (S^1)$:
$$ e=ie^{ikt} d/dt
\tag 4
$$
The commutation relations of the Lie algebra $\operatorname {Vect}_{\Bbb C}
(S^1)$ have the following form
$$
[e_j,e_k]=(j-k)e_{j+k}
\tag 5
$$
in the basis $e_k$.

In 1968 I.M. Gelfand and D.B. Fuchs \cite {7} discovered that $\operatorname
{Vect} (S^1)$ possesses a non-trivial central extension. The corresponding
2-cocycle is
$$
c(u,v)=\int v'(t)du'(t)
\tag 6
$$
or, equivalently,
$$
c(u,v)=\vmatrix
v'(t_0) &u'(t_0)\\
v''(t_0) &u''(t_0)
\endvmatrix. \tag 7
$$
This central extension was independently discovered by M. Virasoro \cite {8}
and named after him. Let us denote the Virasoro algebra by vir. Its
complexification, which is also called the Virasoro algebra, will be denoted
$\operatorname {vir}_{\Bbb C}$. As a vector space vir is generated by the
vectors $e_k$ and the central element $c$. The commutation relations have
the form
$$
[e_j,e_k]=(j-k)e_{j+k}+\delta(j+k)\frac{j^3-j}{12}c.
\tag 8
$$
The infinite-dimensional group Vir corresponding to the
Lie algebra vir is a central extension of the group
$\operatorname {Diff}(S^1)$. The corresponding 2-cocycle was calculated by
R.~ Bott \cite {9}. This cocycle can be written as
$$
c(g_1,g_2)=\int \operatorname {log}(g'_1\circ g_2)\operatorname {log}(g'_2).
\tag 9
$$
The group Vir is
called the Virasoro--Bott group.

There are no groups corresponding to the Lie algebras $\operatorname
{Vect}_{\Bbb C}(S^1)$ or $\operatorname {vir}_{\Bbb C}$, but one can
consider the following construction due to Yu.~A.~Neretin,
M.~L.~Kontzevich, and G. Segal.

Let us denote by ${\operatorname LDiff}^{\Bbb C}_+(S^1)$ the set of all
analytic mappings $g:S^1\mapsto {\Bbb C}\backslash \{0\}$ such that $g(S^1)$ is
a Jordan curve surrounding zero, the
orientations of $S^1$ and $g(S^1)$ are the same, and $g'(e^{i\theta})$ is
everywhere different from zero. ${\operatorname LDiff}^{\Bbb C}_+(S^1)
$ is a local group \cite {10}.

Let $\operatorname{LNer}\subset{\operatorname
LDiff}^{\Bbb C}_+(S^1)$ be the local subsemigroup of mappings $g$ such that
$|g(e^{i\theta})|<1$. As was shown by Yu.~A.~Neretin \cite {10}, the
structure of a local semigroup on LNer extends to the structure of a
global semigroup Ner.

There exist at least two constructions of the semigroup Ner.

The first construction (Yu.~A.~Neretin \cite {10}). An element of Ner is a
formal product $$
p\cdot A(t) \cdot q
\tag 10
$$
where $p,q\in {\operatorname Diff}_+(S^1)$, $p(1)=1$, $t>0$, $A(t):C
\mapsto C$, $A(t)z=e^{-t}z$.

To define the multiplication in Ner one must describe the rule to transform
the formal product $A(s)\cdot p \cdot A(t)$ to the form (10).

{\bf A.} Let $t$ be so small that the diffeomorphism $p$ extends
holomorphically to the annulus $e^{-t} \leqslant |z| \leqslant 1$. Then the
product $g=A(s)pA(t)$ is well-defined. Let $K$ be the domain bounded
by $S^1$ and $g(S^1)$. Let $Q$ be the canonical conformal mapping of $K$ onto
the annulus $e^{-t'} \leqslant |z| \leqslant 1$, normalized by the condition
$Q(1)=1$. Then $g=p'\cdot A(t) \cdot q'$, where $p'=Q^{-1}|_{S^1}$
and $q'$ is determined by the identity
$$
A(s) \cdot p \cdot A(t)=p' \cdot A(t') \cdot q'. \tag 11
$$
{\bf B.} For an arbitrary $t$ there exists a suitable $n$ such that the product
$$
A(s) \cdot p \cdot A(t)=(\dots(A(s) \cdot p \cdot A(t/n))A(t/n)\dots)A(t/n)
\tag 12
$$
can be calculated. It can be shown that the product does not depend on the
choice of the representation (12) and is associative \cite {10}.

The second construction (M.~L.~Kontsevich \cite {11} and G.~Segal \cite {12}).
An element $g$ of the semigroup Ner is a triple $(K,p,q)$, where $K$ is a
Riemann surface with boundary $\partial K$ such that $K$ is biholomorphically
equivalent to an annulus and $p,q:S^1\mapsto \partial K$ are fixed
parametrizations of the components of $ \partial K$. Two elements
$g_i=(K_i,p_i,q_i)$, $i=1,2$, are equivalent if there exists a conformal
mapping
$R:K\mapsto K$ such that $p_2=Rp_1$, and $q_2=Rq_1$. The product of two
elements
$g_1$ and $g_2$ is the element $g_3=(K_3,p_3,q_3)$, where
$$
K_3=K_1 \bigsqcup_{q_1(e^{it})=p_2(e^{it})} K_2,
$$
$p_3=p_1$, and $q_3=q_2$.

This construction admits a slight modification \cite {13}. Let us consider
the semigroup $\overline {\operatorname{Ner}}$ whose elements $g$ are
pairs $(p^+_g, p^-_g)$, where $p^{\pm}_g:D_\pm \mapsto {\Bbb C}$, $D_\pm=\{
z:|z|^{\pm 1}\leqslant 1\}$, such that $p^+_g(D^0_+)\cap
p_g^-(D^0_-)=\varnothing$. Two elements $g_1$ and $g_2$ are equivalent
if there exists a biholomorphic mapping $R:\overline {\Bbb C}\mapsto\overline
{\Bbb C}$ such that $p^+_{g_2}=Rp^+_{g_1}$ and $p^-_{g_2}=Rp^-_{g_1}$. The
product is defined by analogy with the previous construction.

The Neretin semigroup $\bar {\operatorname{Ner}}$ possesses a central
extension.  The corresponding cocycle was calculated by Yu.~A.~Neretin
\cite{13}:
$$
\aligned
c(g_1,g_2)&=\oint\log(p^+_2)'(z)d\log\frac{(p^+_3)'(g_1(z))}{(p^+_1)'(g_1(z))}-
\oint\frac{\log(p^+_2)(z)}{z}d\log\frac{p^+_3(g_1(z))}{p^+_1(g_1(z))}\\
&\qquad +
\oint\log(p^-_2)(g_2(z))d\log\frac{p^-_3(z)}{p^-_1(z)}
-\oint\log\frac{p^-_2(g_2(z))}{g_2(z)}d\log\frac{p^-_3(z)}{p^-_1(z)}.
\endaligned
\tag 13
$$
\subhead 1.2. The flag manifold of the Virasoro-Bott group
\endsubhead
The flag manifold $M$ of the Virasoro-Bott group is a homogeneous space with
transformation group $\operatorname {Diff}_+(S^1)$ and isotropy group $S^1$.
There exist several different realizations of this manifold \cite{14--18}.

{\it Algebraic realization\/}. The space $M$ can be realized as a conjugacy
class in the group $\operatorname {Diff}_+(S^1)$ or in the Virasoro-Bott group
Vir \cite{15}. It should be mentioned that $M$ can also be realized as
the quotient $\operatorname{Ner}/{\operatorname{Ner}}^{\circ}$ of the Neretin
semigroup by the subsemigroup ${\operatorname{Ner}}^{\circ}$ of
elements admiting
holomorphic extension to $D_-$.

{\it Probabilistic realization\/}. Let $P$ be
the space of real probability measures $\mu = u(t)\,dt$
with smooth positive density $u(t)$ on $S^1$. The group
$\operatorname {Diff}_+(S^1)$
naturally acts on $P$ by the formula
$$
g\: u(t)\,dt \mapsto u(g^{-1}(t))\,dg^{-1}(t). \tag 14
$$
The action is transitive and the stabilizer of the point $(2\pi)^{-1}\,dt$ is
isomorphic to $S^1$. Hence,
$P$ can be identified with $M$.

{\it Orbital realization\/}. The space $M$ can be considered as an orbit of the
coadjoint representation of $\operatorname {Diff}_+(S^1)$ or Vir \cite{17,
18}. Namely, the elements of the dual space $\operatorname{vir}^*$ of the
Virasoro algebra vir can be identified with the pairs $(p(t)\,dt^2,b)$; the
coadjoint action of Vir has the form
$$
K(g)(p,b)=(gp-bS(g),b), \tag 15
$$
where
$$
S(g)=\frac{g'''}{g'}-\frac32 \Big(\frac{g''}{g'}\Big)^2
$$
is the Schwarzian (the Schwarz derivative). The orbit of the point $(a\cdot
dt^2,b)$  coincides with $M$ provided that $a/b \neq -n^2/2$, $n=1,2,3,\dots$.
Therefore, a family $\omega_{a,b}$ of symplectic structures  is defined on $M$.

{\it Analytic realization\/}. Let us consider the space $S$ of univalent
functions on the unit disk $D_+$ \cite {19--21}. The Taylor coefficients
$c_1,c_2,c_3,c_4,\dots$ in the expansion
$$
f(z)=z+c_1z^2+c_2z^3+c_3z^4+\dots+c_nz^{n+1}\dots \tag 16
$$
form a coordinate system on $S$. It was shown in \cite {14} that $S$
can be naturally identified with $M$. In the
coordinate system $\{c_k\}$ the action of the Lie algebra
$\operatorname {Vect}_{\Bbb C}(S^1)$ on $M$ has the following form:
$$
{\frak{L}}_v(f(z))=if^2(z)\oint\Big(\frac{wf'(w)}{f(w)}\Big)^2
\frac{v(w)}{f(w)-f(z)}\frac {dw}{w} \tag 17A
$$
or
$$
\aligned
L_p&=\frac{a}{ac_p}+\sum_{k\geqslant 1}(k+1)c_k\frac a{ac_{k+p}},\quad p>0\\
L_0&=\sum_{k\geqslant 1}kc_k \frac a{ac_{k}}\\
L_{-1}&=\sum_{k\geqslant 1}((k+2)c_{k+1})\frac a{ac_{k}}\\
L_{-2}&=\sum_{k\geqslant
1}((k+3)c_{k+2})-(4c^2-c^2_1)c_k-B_k) \frac a{ac_{k}}\\
L_{-p}&=\frac{(-1)^p}{(p-2)!}ad^{p-2}(L_{-1})L_{-2}
\endaligned
\tag 17B
$$
where $B_k$ are the Laurent coefficients of the
function $1/(wf(w))$. The symplectic structure $\omega_{a,b}$ together
with the complex structure on $M$ determines a K\"ahler metric $w_{a,b}$.
More detailed information can be found in \cite {3, 15,
16, 22, 23}.

It should be mentioned that the space $M$ can be realized as a space of complex
structures on loop manifolds \cite{24}.
\subhead 1.3. Non-Euclidean geometry of mirrors
\endsubhead
Points and Lagrangian submanifolds are the basic
elements of symplectic  geometry \cite {25-28}. However the space of all
Lagrangian submanifolds is infinite-dimensional and hard to visualize, which is
not convenient. The K\"ahler geometry on the flag manifold $M$ of the
Virasoro-Bott group permits us to select a handy subset in the set of
all Lagrangian submanifolds \cite {29}.

By a K\"ahler subsymmetry \cite {29, 4--6} we mean an involutory
antiautomorphism
of $M$. The set of fixed points of a subsymmetry (a {\it mirror\/}) is a
completely geodesic Lagrangian submanifold. Points and mirrors are the basic
elements of the geometry to be described. The set of all mirrors forms a
symmetric space, the non-Euclidean Lagrangian (Lagrange Grassmannian) $\Lambda
(M)$, independent of a choice of the non-Einsteinian K\"ahler
$\operatorname {Diff}_+ (S^1)$-invariant metric $w_{a,b}$ $(b/a\neq-13)$ on
the space $M$ \cite{24, 3, 29}.

Consider the probabilistic realization of $M$ and the set $A(M)$  of
measures of the form $\delta_a(t)dt$. The set $A(M)$ is called the
absolute \cite{3, 29}. The absolute $A(M)$ is isomorphic to $S^1$. Let us
introduce the parallelism relation  on the non-Euclidean Lagrangian. Two
mirrors are said to be parallel if they pass through the same point of
the absolute. The following analog of the Lobachevskii axiom holds:
for any point of $M$ and any point of
$A(M)$ there exists exactly one mirror passing through these points.

Let us also introduce the non-Euclidean Lagrangian $\Lambda_+(M)$
\cite{29}. The elements of $\Lambda_+(M)$ are the oriented mirrors $V_a$.
An oriented mirror is a pair $(V,a)$, where $V$ is a mirror
and $a\in V$ a point of the absolute.

It should also be mentioned that the Maslov index $m(U_a,V_b,W_c)$
\cite{30} of three oriented mirrors $U_a,V_b$, and $W_c$ is equal to $1$ if the
orientations of the triple $(a,b,c)$ and
of the circle $S^1\simeq A(M)$ are the same and to $-1$ if these orientations
are opposite.

\head 2. The infinite-dimensional geometry of the skeleton of the flag manifold
of the Virasoro-Bott group\endhead

\subhead 2.1. An equivariant mapping  of the flag manifold of the
Virasoro-Bott group into the infinite-dimensional classical domain of third
type\endsubhead
Let $H_0$ be the space of smooth real-valued
$1$-forms $u(\exp(it))\,dt$ on the
circle such that
$$
\int u(\exp(it))\,dt=0.
$$
Let $H_0^{\Bbb C}$ be its complexification, and let $(H_0^{\Bbb C})_+$ and
$(H_0^{\Bbb C})_-$ be the transversal spaces
consisting of $1$-forms possessing holomorphic extensions to the disks $D_+$
and $D_-$, respectively. Let $\Cal O(S^1)$, $\Cal O(D_+)$, and $\Cal O(D_-)$
denote the spaces of holomorphic functions on $S^1$, $D_+$, and $D_-$,
respectively. Then the space $H_0^{\Bbb C}$ is isomorphic to the quotient of
$\Cal O(S^1)$ by constants:
$$
f(z)\in \Cal O(S^1)\mapsto df(z)\in H_0^{\Bbb C}.
$$
Under this isomorphism the spaces $(H_0^{\Bbb C})_+$ and $(H_0^{\Bbb C})_-$ are
identified with the quotients of $\Cal O(D_+)$ and $\Cal O(D_-)$, respectively,
by constants. Consider the completion $H^{\Bbb C}$ of the space
$H_0^{\Bbb C}$ with respect to the norm
$$
\|u\|=\sum_n|u_n|^2/n,\quad
u(z)=\sum_n u_n z^n,
$$
and let $H_+^{\Bbb C}$ and $H_-^{\Bbb C}$ be the corresponding completions of
$(H_0^{\Bbb C})_+$ and $(H_0^{\Bbb C})_-$. The space $H^{\Bbb C}$ is equipped
with the symplectic and pseudo-Hermitian structures defined by
$$
\aligned
(f(z),g(z))&=\oint f(z)\,dg(z), \\
\langle f(z),g(z)\rangle&=\oint f(z)\,\overline{dg(z)}, \quad f,g\in\Cal
O(S^1).
\endaligned
\tag 18
$$
Let us denote the invariance groups of these structures by
$\operatorname{Sp}(H^{\Bbb C},\Bbb C)$ and $U(H^{\Bbb C}_+,H^{\Bbb
C}_-)$, respec\-tively; also, let $\operatorname{Sp}(H,\Bbb R) =
\operatorname{Sp}(H^{\Bbb
C},\Bbb C)\cap U(H_+,H_-)$.

Consider the Grassmannian $\operatorname{Gr}(H^{\Bbb C})$, that is, the set of
all complex Lagrangian subspaces in $H^{\Bbb C}$ \cite{48, 42}.
The space $\operatorname{Gr}(H^{\Bbb C})$ is an infinite-dimensional
homogeneous
space with transformation group $\operatorname{Sp}(H^{\Bbb
C},\Bbb C)$. Consider the action of the subgroup $\operatorname{Sp}(H,\Bbb R)$
on
$\operatorname{Gr}(H^{\Bbb C})$. The orbit of the point $H^{\Bbb C}_-$ is an
open (in a suitable topology) subspace $\Cal R$ in $\operatorname{Gr}(H^{\Bbb
C})$ isomorphic to $\operatorname{Sp}(H,\Bbb R)/U$, where $U$ is the group of
operators on $H^{\Bbb C}=H^{\Bbb C}_+\oplus H^{\Bbb C}_-$ of the form
$A\oplus \widetilde A$, $A\in U(H^{\Bbb C}_+)$, $\widetilde A\in U(H^{\Bbb
C}_-)$; here the mapping $A\mapsto \widetilde A$ from $U(H^{\Bbb C}_+)$ to
$U(H^{\Bbb C}_-)$ is induced by the inversion $z\mapsto 1/z$.

The manifold $\Cal R$ is an infinite-dimensional classical homogeneous domain
of
third type \cite{31}. The manifold $\Cal R$ is mapped in the linear space
$\operatorname{Hom}(H^{\Bbb C}_-,H^{\Bbb C}_+)$ in such a way that the
elements of $\Cal R$ are represented by symmetric matrices $Z$ satisfying
$E-Z\bar Z>0$.

The mapping of $M$ into $\Cal R$ is described in
\cite{16, 22, 23, 3}. Namely, the representation of
$\operatorname{Diff}_+(S^1)$ in $H$ defines a monomorphism
$\operatorname{Diff}_+(S^1) \mapsto \operatorname{Sp}(H,\Bbb R)$. Hence,
$\operatorname{Diff}_+(S^1)$ acts on $\Cal R$. The orbit of the initial point
under this action coincides with
$\operatorname{Diff}_+(S^1)/\operatorname{PSl}(2,\Bbb R)$. Therefore, we have a
mapping $M\mapsto \Cal R$. The explicit form of this mapping can be found in
\cite{16, 23} (see \cite{3}). The matrix $Z_f$ corresponding to a univalent
function $f\in S$ is called the Grunskii matrix, and the mapping $S\simeq
M \mapsto \Cal R\hookrightarrow
\operatorname{Gr}(H^{\Bbb C})$ is the Krichever mapping \cite{41,
42}.

It is well known \cite{32} that the skeleton of a finite-dimensional classical
domain of third type consists of all symmetric unitary matrices. Thus we
can regard the set of all symmetric unitary operators from $H$ to $H$ as the
skeleton of $\Cal R$. By the skeleton of the space $S$ of univalent
functions is defined as the set of functions whose Grunskii matrices are
unitary.
Accordingly to Milin's theorem \cite{33, 34}, the skeleton of $S$
consists of all univalent functions $f$ such that
$\operatorname{mes}(\Bbb C\setminus f(D_+^0))=0$. Let us investigate the
structure
of $S$ more systematically.

\subhead 2.2. The geometry of the skeleton of the space $S$\endsubhead
Consider the $\Bbb R$-analytic space $E$ whose elements are cuts of the complex
plane $\Bbb C$ with one end at infinity such that the conformal radius of
$\Bbb C\setminus K$ with respect to zero is equal to one. Consider the
mapping $E\to \Lambda_+(M)$ defined as
$$
f(z)\mapsto (s,a), \tag 19
$$
where
$$
\aligned
f(D_+^0)\sqcup K &= \bar C, \quad f(0)=0, \quad f'(0)=1,\\
f(a) &=\infty,\quad f(s(z))=f(z).
\endaligned
\tag$19'$
$$

\proclaim{Theorem 1A} $E\simeq \Lambda_+(M)$.
\endproclaim
\demo{Proof} Note that $\Lambda_+(M) = \{(s,a)\:a\in
\operatorname{Diff}_-(S^1)$, $s^2=id$, $a\in S^1$, $s(a)=a$\}. It is clear
that the mapping (19) is an embedding. Let us prove that it is a surjection.
To this end consider an arbitrary element $(s,a)$ of $\Lambda_+(M)$ and
construct the manifold $D_+^s = D_+/(z=s(z))$. Then $D_+^s$ is topologically
equivalent to the Riemann sphere, and therefore $D_+^s$ and $\overline C$ are
equivalent as complex manifolds. Hence, there exists unique mapping
$f\:D_+^s\mapsto \Bbb C$ such that $f(0)=0$, $f'(0)=1$, and $f(a)=\infty$. The
composition of $f$ with the natural projection $D_+\to D_+^s$ is a
function representing the element of $E$ corresponding to the pair $(s,a)$.
\enddemo

Let us now embed $E$ in the skeleton of $S$. Note that the group
$\operatorname{Diff}(S^1)$ naturally acts on the skeleton. On the other hand,
Theorem 1A defines the structure of a $\operatorname{Diff}_+(S^1)$-homogeneous
space on $E$. The question is whether the embedding of $E$ in the skeleton of
$S$ is $\operatorname{Diff}_+(S^1)$-equivariant.

\proclaim{Theorem 1B} The embedding of $E\simeq\Lambda_+(M)$ in the skeleton of
$S$ is $\operatorname{Diff}_+(S^1)$-equivariant. \endproclaim
\demo{Proof} Note that the action of $\operatorname{Diff}_+(S^1)$ on the
skeleton of $S$ can be reduced to the subspace $E$. Namely, formulas (17A)
that define the infinitesimal action correspond to analytic variations of cuts
in $E$.  One obtains the action of
$\operatorname{Diff}_+(S^1)$ on $E$ by exponentiating these deformations.

It is necessary to verify that this action is the same as defined in
the theorem.

Note that the group $\operatorname{Diff}(S^1)$ naturally
acts on $\operatorname{Gr}(H^{\Bbb C})$ preserving $R$ and its skeleton. On
the other hand, $\operatorname{Diff}(S^1)$ preserves $S\simeq
\operatorname{Diff}_+(S^1)/S^1$ and therefore preserves its skeleton. This
statement is true for each K\"ahler subsymmetry of $S$ (which is an involutory
element of $\operatorname{Diff}_-(S^1)$). Consider an arbitrary
subsymmetry $s\in\Lambda(M)$. Its mirror $V_s$ extends analytically to be an
element of $\operatorname{Gr}(H^{\Bbb C})$. The set $V_s\cap E$ consists of
exactly two points that correspond to different elements $(s,a)$ and
$(s,b)$ of $\Lambda_+(M)$ under the isomorphism (19). Thus, the cited actions
of
$\operatorname{Diff}_+(S^1)$ on $E$ are the same.
\enddemo

Note that $\Lambda_+(M)$ is a symmetric space,
$$
(s_1,a_1)(s_2,a_2) = (s_1s_2s_1, s_1(a_2)),
$$
with trasvection group $\operatorname{Diff}_+(S^1)$ and isotropy group $G_0 =
\{g\in\operatorname{Diff}(S^1)\: g(1)=1,\,\,\overline{g(z)}=g(\bar z)\}$.
The tangent space $V$ to $\Lambda_+(M)$ at the point $(s_-,1)$,
$s_-(z)=\bar z$, can be identified with the space of odd vector fields on
$S^1$. The directions in $V$ invariant with respect to $G_0$ and
determined by the generalized   vectors $\delta_1(t)d/dt$ and
$\delta_{-1}(t)d/dt$ give rise to nonholonomic generalized invariant direction
fields $\xi_+$ and $\xi_-$ on $\Lambda_+(M)$. Let $\Cal O(E)$ and
$\Cal O(\Lambda_+(M))$ be the structure rings of $E$ and $\Lambda_+(M)$,
respectively; we have $\Cal O(\Lambda_+(M)/\xi_-)= \{f\in
\Cal O(\Lambda_+(M))\:\xi_-f=0\}$.
\proclaim{Theorem 1C}
$(E,\Cal O(E))\simeq (\Lambda_+(M),\Cal O(\Lambda_+(M)/\xi_-))$.
\endproclaim
\demo{Proof} As was mentioned, the
tangent space to $\Lambda_+(M)$ at the point $(s_-,1)$, $s_-(z)=\bar z$, can be
identified with the space of odd vector fields on $S^1$. Under the isomorphism
(19) the point $(s_-,1)\in \Lambda_+(M)$ corresponds to the Koebe function
\cite{19, 2}
$$
k(z) = z/(1-z)^2.
$$
The Killing fields on $E$ defined by (17A) and vanishing at the point
$k(z)$ are just the odd vector fields and $\xi_-$. Since $E$ and $\Lambda_+(M)$
are homogeneous spaces, it follows that $\Cal O(E) = \Cal
O(\Lambda_+(M)/\xi_-)$.
\enddemo

Let us now consider the mapping $M\mapsto\Gamma_{cl}(\Lambda_+(M))$, where
$\Gamma_{cl}(\Lambda_+(M))$ is the space of all closed geodesics on
$\Lambda_+(M)$, that assigns to each point $x\in M$ the set of all oriented
mirrors passing through $x$. This mapping is clearly an
isomorphism. Namely, consider a closed geodesic on $\Lambda_+(M)$. This
geodesic
is a symmetric space with transvection group isomorphic to some subgroup
$S^1\subset \operatorname{Diff}_+(S^1)$. This subgroup is the stabilizer  of
a suitable point of the flag manifold $M$.

Under the identification of $M$ and $\Gamma_{cl}(\Lambda_+(M))$ the
symplectic structure on $M$ has the form
$$
\omega_x(X,Y) = \int_{\gamma_x}(AX,Y)\,d\tau, \tag 20
$$
where $\tau$ is the natural parameter on $\gamma_x$ and $X,Y$ are Jacobi
fields orthogonal to the field $\dot{\tau}$ in the unique (up to
a real factor) invariant degenerate pseudo-Riemannian metric on
$\Lambda_+(M)$; $A=a\nabla + b\nabla^3$, where $\nabla$ is the covariant
derivative along $\gamma_x$ (cf. \cite{35}).

Consider the class $\Cal O^0(S)$  of holomorphic functionals on $S$ that
admit analytic extension to $E$.
\proclaim{Theorem 2} $\Re(\Cal O^0(S))\simeq \Cal O(\Lambda_+(M)/\xi_-)$. The
isomorphism is given by the Poisson-type formula
$$
F(f) = \int\limits_{f_\tau\in\gamma_x} F(f_\tau)\,d\tau,
$$
where $f$ is univalent function of class $S$ corresponding to the
point $x$ of the homogeneous manifold $M$,
$\gamma_x$ is a geodesic on $\Lambda_+(M)=E$, $f_\tau$ are the univalent
functions from $E$ corresponding to the points of the geodesic $\gamma_x$,
$\tau$ is the natural parameter on $\gamma_x$, and $F$ is a functional from
$\Re
\Cal O^0(S)$.
\endproclaim
\demo{Proof} Let $f$ be the function of class $S$ of the form $f(z)\equiv
z$; then the univalent functions $f_\tau$ are the Koebe functions
$k_\tau(z)=z/(1-e^{i\tau}z)$  \cite{19, 20}. Consider the
analytic disk $D=\{f_a(z) = z/(1-az)^2,\,\,a\in D_+\}$ lying in $S$. The
restriction $F$ of an analytic functional from $\Re
\Cal O^0(S)$ to the disk $D$ is a harmonic function on $D$, and, therefore,
$$
F(f_0)=\int_{f_\tau}F(f_\tau)\,d\tau
$$
for all $f\in S$.
\enddemo

\head 3. The  infinite-dimensional geometry of the universal deformation of the
complex disk
\endhead
\subhead
3.1. The universal deformation of the complex disk
\endsubhead
The space $A$ of the universal deformation of the complex disk is the
set of pairs
$$
(f,w),\quad f\in S,\quad w\in C\setminus ((f(D_+)^{-1}\cup\{0\}), \tag 21A
$$
and the projection onto the base has the form
$$
(f,w)\mapsto f. \tag 21B
$$
The natural action of the Virasoro algebra $\operatorname{vir}_{\Bbb C}$ on
$A$,
introduced in \cite{37,  38} and explicated in \cite{36}, has the form
$$
{\frak L}_v(f(z), a) = i\oint\Big(\frac{wf'(w)}{f(w)}\Big)^2\frac{v(w)\,dw}{w}
\Big(\frac{f^2(z)}{f(w)-f(z)}, \frac{a}{1-af(w)}\Big) \tag 22A
$$
or
$$
\align
L_p &= \frac{\partial}{\partial c_p}   + \sum\limits_{k\geqslant 1} (k+1)c_k
\frac{\partial}{\partial c_{k+p}},\quad p>0\\
L_0 &=\sum\limits_{k\geqslant 1} kc_k
\frac{\partial}{{\partial c}_k} + w\partial_w\\
L_{-1} &=\sum\limits_{k\geqslant 1} ((k+2)c_{k+1} -
2c_1c_k)\frac{\partial}{\partial c_k} + 2c_1w \frac{\partial}{\partial w} +
w^2\frac{\partial}{\partial w} + w^2\frac{\partial}{\partial w}\\
L_{-2} &=\sum\limits_{k\geqslant 1} ((k+3)c_{k+2} - (4c_2-c_1^2)c_k - B_k)
\frac{\partial}{\partial c_k} + (4c_2-c_1^2)w\frac{\partial}{\partial w} +
3c_1w^2\frac{\partial}{\partial w} + w^3\frac{\partial}{\partial w}\\
L_{-p} &= \frac{(-1)^p}{(p-2)!} ad^{p-2}(L_{-1})L_{-2},\quad p\geqslant 3.
\endalign
$$

This action can be exponentiated to yield an action of the Neretin semigroup.
Thus the universal deformation is identified with the quotient
$\operatorname{Ner}/\operatorname{Ner}^{\circ\circ}$, where
$\operatorname{Ner}^{\circ\circ}$
is the subsemigroup of codimension $1$ in $\operatorname{Ner}^{\circ}$
consisting
of mappings $g\in \Cal O(D_-)$ with some prescribed fixed point.

The action of $\operatorname{Diff}_+(S^1)$ on the
base $M$ can be lifted to the  universal deformation space $A$.

\subhead
3.2. The universal deformation space of the complex disk as a
subsymmetric space \endsubhead
\definition{Definition 1} \cite{4, 5} A pair $(\Cal X,\Sigma)$, where $\Cal X$
is a space and $\Sigma$ the set of its involutive automorphisms
(antiautomorphisms) (involutions), is called a {\it subsymmetric space\/} if
a mapping
$$
x\mapsto s_x
$$
from $X$ to $\Sigma$ is given such that

a) $s_x(x) = x$, $s_xs_ys_x=s_{s_xy}$,

b) if for some $s\in\Sigma$ $s(x)=x$, then $s=s_x$.
\enddefinition
\proclaim{Lemma} Each subsymmetry of $M$ can be extended to a subsymmetry of
$A$,
and this extension is $\operatorname{Diff}_+(S^1)$-equivariant.
\endproclaim
\demo{Proof} The assertion of the lemma follows from the fact that the
elements of $\operatorname{Diff}_+(S^1)$ can be regarded as
automorphisms of the semigroup $\operatorname{Ner}$. This fact follows from
the $\operatorname{Diff}_+(S^1)$-invariance of the tangent cone ner to the
semigroup Ner, which lies  in the Lie algebra $\operatorname{Vect}_{\Bbb
C}(S^1)$.
\enddemo

\proclaim{Theorem 3A} The pair $(A,\Lambda(M))$ is a subsymmetric space.
\endproclaim
\remark{Remark} The universal deformation space of the complex disk
is projected on the symmetric space $\Lambda(M)$. Moreover, since the
subsymmetry mirrors on $A$ consist of  two connected
components, the universal deformation space of the complex disk is
projected on the symmetric space $\Lambda_+(M)$.
\endremark

In view of this fact the non-Euclidean Lagrangian $\Lambda_+(M)$
can be viewed as
a ``the universal deformation space of the circle."

\proclaim{Theorem 3B} $\Re(\Cal O(A))\simeq \Cal O(\Lambda_+(M))$.
\endproclaim
\demo{Proof} The mapping of $\Re(\Cal O(A))$ into $\Cal O(\Lambda_+(M))$ has
the
form
$$
F(s,a) = \lim_{f_n\to f,\,\, w\to b} F(f_n, b_w),
$$
where $f$ is the function corresponding to the cut in $E$ determined by the
element $(s,a)\in\Lambda_+(M)$, and $1/b$ is the second end of the cut
$f(S^1)$.
The surjectivity of this mapping follows from Theorem 2 and formulas
(22A).
\enddemo

Note that $A$ is a symplectic \cite{39, 40} and, moreover, a
K\"ahler manifold. The K\"ahler structure can be defined by Bergman's kernel
function, which is the exponential of K\"ahler's potential.

This  kernel function is the product of the lift of Bergman's
kernel function on the base $M$ and fiberwise Bergman's
kernel function. The subsymmetries of $A$ are K\"ahler antiautomorphisms.
In particular, the element $s_-$ defines the subsymmetry of $A$ of the form
$$
(f(z),w)\mapsto (\overline{f(\bar z)}, \bar w).
$$
Mirrors of subsymmetries are Lagrangian submanifolds, and the projection of $A$
to $\Lambda_+(M)$ defines a real $\operatorname{Diff}_+(S^1)$-invariant
polarization (see \cite{43}) on $A$.

The author is grateful to M.~A.~Semenov-Tyan-Shanskii, A.~S.~Schwarz,
A.~S.~Fedenko, L.~V.~Sa\-binin, I.~M.~Milin, G.~I.~Ol$'$shanskii,
A.~M.~Perelomov, A.~Yu.~Morozov, M.~L.~Kontsevich, Yu.~A.~Neretin,
P.~O.~Mikheev,
E.~G.~Emel$'$yanov, A.~O.~Ra\-dul, A.~A.~Roslyi, A.~Yu.~Alekseev, and
B.~A.~Khesin
for useful discussions of some aspects of the paper. The author thanks
A.~N.~Rudakov and the participants of his seminar for the attention. The author
is indebted to V.~N.~Kolo\-kol$'$tsov for fruitful discussions.

\enddocument